\documentclass{amsart}
\usepackage{graphicx,amssymb,epsfig}
\newtheorem{theorem}{Theorem}[section]
\newtheorem{proposition}[theorem]{Proposition}
\newtheorem{lemma}[theorem]{Lemma}
\newtheorem{corollary}[theorem]{Corollary}
\newtheorem{remark}[theorem]{Remark}
\newtheorem{example}[theorem]{Example}

\newtheorem{definition}[theorem]{Definition}

\newcommand{\bth}{\begin{theorem}}
\newcommand{\bpr}{\begin{proposition}}
\newcommand{\epr}{\end{proposition}}
\newcommand{\bco}{\begin{corollary}}
\newcommand{\eco}{\end{corollary}}
\newcommand{\bde}{\begin{definition}}
\newcommand{\ede}{\end{definition}}
\newcommand{\ble}{\begin{lemma}}
\newcommand{\ele}{\end{lemma}}
\newcommand{\bre}{\begin{remark}\rm}
\newcommand{\ere}{\end{remark}}
\newcommand{\bex}{\begin{example}\rm}
\newcommand{\eex}{\end{example}}

\def\la#1{\hbox to #1pc{\leftarrowfill}}
\def\ra#1{\hbox to #1pc{\rightarrowfill}}

\def\fract#1#2{\raise4pt\hbox{$ #1 \atop #2 $}}

\def\lrar{{\ra 2}}

\newcommand{\dis}{\displaystyle}
\def\bbc{\mathbb C}
\def\bbp{\mathbb P}

\def\bbz{\mathbb Z}
\def\bbr{\mathbb R}

\begin{document}

\title{Invariant Spin Structures on Riemann Surfaces}

\author{Sadok Kallel}
\author{Denis Sjerve\ *}
\thanks{*Research partially supported by NSERC
Discovery grant A7218}

\begin{abstract} We investigate the action of the automorphism group of a
  closed Riemann surface on its set of theta characteristics (or spin
  structures).  We give criteria for when an automorphism fixes all spin
  structures, or when it fixes just one. The case of hyperelliptic curves and
  of the Klein quartic are discussed in detail.
\end{abstract}
\maketitle
%******************************************************************************

\section{Introduction}\label{intro}

Spin structures on Riemann surfaces or ``theta characteristics" are classical
objects of great use and interest in all of geometry, topology and physics.
Let $C$ be a closed Riemann surface and write $U(C)$ for its unit tangent
bundle.  A spin structure on $C$ is a cohomology class in $H^1(U(C );\bbz_2)$
whose restriction to each fiber is a generator of $H^1(S^1;\bbz_2)$.  There is
of course a great number of ways to see what a spin structure is, some of
which we will encounter in the course of this paper.

Every Riemann surface $C$ has $2^{2g}$ distinct spin structures where $g$ is
the genus of $C$.  Generally the set of spin structures $Spin(C )$ is an
affine space associated with $H^1(C,\bbz_2)$ and hence there is a
(non-canonical) $1-1$ correspondence between $Spin(C )$ and $H^1(C;\bbz_2)$.

The group of automorphisms $Aut(C)$ acts on $H^1(U(C);\mathbb Z_2)$ and hence
on $Spin(C )$ by pullback. This action and the action on $H^1(C;\bbz_2)$ are
however not compatible under the above $1-1$ correspondence, but are closely
related.  They differ by a term involving the intersection pairing.

In this paper we study invariant spin structures on compact Riemann surfaces.
Motivation for this work comes from a question by Jack Morava to the authors
and from a paper of Atiyah \cite{atiyah}.  In that paper Atiyah proves that if
$f:C\to C$ is an automorphism of a compact Riemann surface $C$ then there
necessarily exists a spin structure $\xi$ that is invariant under $f:C\to C.$

This raises natural questions.

\begin{enumerate}
\item How many spin structures are invariant under  $f:C\to C?$
\item How large can the isotropy subgroup of a given
  spin structure in $Aut(C)$ be?
\item Is there an  automorphism $f:C\to C$ that leaves only one spin structure invariant?
\item Is there a non-trivial automorphism $f:C\to C$ that leaves every spin
  structure invariant?
\end{enumerate}

We will give complete answers to these questions and more. 
The following is an answer to (4) and shows how restrictive that question is.

\begin{theorem}\label{mainth3}
A non-trivial automorphism $f:C\to C$ leaves every spin structure invariant if,
and only if, $C$ is hyperelliptic and $f$ is the hyperelliptic involution.
\end{theorem}

In fact one can determine precisely the number of invariant spin structures
for an automorphism $f$ in terms of its symplectic action on homology.  Let
$A\in SL_{2g}(\mathbb Z)$ denote the $2g\times 2g$ matrix representing the
induced isomorphism $f_*:H_1(C;\mathbb Z)\to H_1(C;\mathbb Z),$ with respect
to some basis, and let $\bar{A}$ is its mod 2 reduction. If $h$ is the
dimension of the eigenspace of $\bar{A}$ associated to the eigenvalue $1$,
then the number of $f$-invariant spin structures is $2^h$. This is our answer
to (1).

Note that if $f$  has order $n$, then necessarily $h\geq 2k$ where $k$ is 
the genus of the orbit surface $C/\mathbb Z_n.$ 
This is because $2k$ is the dimension
of the eigenspace of $A$ associated to the eigenvalue $1.$

In contrast with theorem \ref{mainth3}, by analyzing the similarity class of
$A$ in $SL_{2g}(\bbz )$ we obtain

\begin{theorem}\label{mainth2}
  Suppose $f:C\to C$ is an automorphism of order $n,$ where $n$ is odd. Then
  $f$ leaves only one spin structure invariant if, and only if, the associated
  orbit surface $C/\mathbb Z_n$ has genus zero.
 \end{theorem}

 The above theorems are based on a short calculation that we give in section
 \ref{spin} using constructions of Johnson \cite{johnson}. Let $A,\bar A$ be
 as above, and let $X=(x_1,\ldots,x_{2g})^T$ be a column vector representing
 an element of $H^1(C;\mathbb Z)$, and let $\bar{X}$ be its mod 2 reduction.
 Then

 \bpr\label{mainth1} The spin structures left invariant by an automorphism
 $f:C\to C$ are in 1-1 correspondence with solutions $\bar{X}$ of the matrix
 equation $(\bar{A}^T-I)\bar{X}^T=0.$ \epr In particular we see that an
 automorphism $f$ of $C$, or more generally an orientation preserving
 diffeomorphism, fixes all spin structures if and only if $f$ acts trivially
 on $H_1(C;\bbz_2)$. This is one form of the main theorem A of Sipe
 \cite{sipe}.

 Another immediate corollary of proposition \ref{mainth1} and of computations
 of \cite{rauch} is that the Klein quartic curve $\mathcal K$ of genus $3$ and
 maximal group of automorphisms has a unique invariant spin structure.
 Describing the geometry behind this structure requires the description of
 spin structures in terms of divisor classes (cf. \S\ref{appendix}, Appendix).

Consider the points $a=[1:0:0], b=[0:1:0]$ and $c=[0:0:1]$ belonging to
$\mathcal K$, where $\mathcal K$ is described as the locus $x^3y+y^3z+z^3x=0$
in $\bbp^2$. Let $K$ be the canonical divisor on ${\mathcal K}$ which is
determined by the intersection of the curve with any hyperplane
$\bbp^1\subset\bbp^2$. In \S\ref{kleincurve} we prove:

\bth\label{mainth4} The divisor class $\theta := 2a+2b+2c - K$ is the unique
spin structure on the Klein curve ${\mathcal K}$ fixed by the entire group
$Aut( {\mathcal K}) = PSL(2,\bold F_7)$.
\end{theorem}

The existence of a unique spin structure on Klein's curve is well-known to
algebraic geometers \cite{dolgachev2}. It comes as follows.  Let $X(p), p>5$
prime, denote the modular curve defined as a compactification of the upper
half-plane by the action of the principal congruence subgroup of level $p$.
The group $G={\rm PSL}(2,\bold F_p)$ acts as a group of automorphisms on
$X(p)$.  Then \cite{adler, dolgachev} the group ${\rm Pic}(X(p))^G$ of
invariant line bundles is infinite with cyclic generator a $(2p-12)$-th root
of the canonical bundle of degree $(p^2-1)/24$.  When $p=7$, $X(7)=K$ and the
generator is a square root of the canonical bundle hence the spin structure in
question.

Note that if $Aut(C,L)$ is the subgroup of all automorphisms of $Aut(C)$ that
fix a spin structure $L$, then the above corollary shows it is possible that
$Aut(C,L)$ reaches the maximum hurwitz bound of $84(g-1)$. This answers our
question (2).

Even though the above basic results do not seem to exist in the literature, we
make a cautious claim to originality since some of what we prove might be
known to algebraic geometers, coated in some orthogonal language to ours.
Complementary results of a similar nature as ours can be found
in \cite{sipe}. A relevant discussion from the point of view of mapping class
groups is in \cite{masbaum}.  Finally a growing number of references on this
question have to do with moduli spaces of spin curves.

This note is organized as follows.  In \S\ref{spin} we collect the background
material we need and prove our main results.  In \S \ref{hyperelliptics} we
give a complete count of invariant spin structures for automorphisms of
hyperelliptic surfaces based on a combinatorial definition of spin structures
due to Mumford. See propositions \ref{hyperinvodd} and \ref{hyperinveven}. In the
last section \S\ref{kleincurve} we discuss the case of quartics in general and
the Klein curve in particular. An appendix clarifies the relationship between
the two different definitions of spin structures we use in this paper.

%******************************************************************************

\section{Spin Structures and Automorphisms}\label{spin}

We assume $C$ is smooth, closed, connected and orientable.  Let
$S^1\stackrel{i}{\to} U(C)\stackrel{\pi}{\to} C$ be the unit tangent bundle.
We adopt Johnson's definition of a spin structure \cite{johnson}; namely this
is a cohomology class $\xi\in H^1(U(C);\bbz_2)$ that restricts to a generator
of $H^1(S^1;\bbz_2)\approx\ \mathbb Z_2$ for every fiber $S^1.$

There are  short exact sequences
\begin{center}$\begin{array}{ccccccccc}\displaystyle
&0&\to & H_1(S^1;\bbz_2)\stackrel{i_*}{\to}& H_1(U(C);\bbz_2)&\stackrel{\pi_*}{\to}&
H_1(C;\bbz_2)&\to&0\\
&0&\leftarrow & H^1(S^1;\bbz_2)\stackrel{i^*}{\leftarrow}& H^1(U(C);\bbz_2)&
\stackrel{\pi^*}
{\leftarrow}&H^1(C;\bbz_2)&\leftarrow&0
\end{array}$
\end{center}
from which one deduces that the set $Spin(C)$ of spin structures on $C$ is in
$1-1$ correspondence with $H^1(C;\bbz_2) \approx \bbz_2^{2g},$ and therefore
there are $2^{2g}$ spin structures on $C.$ This correspondence is not a group
isomorphism and $Spin(C)$ is in fact the non-trivial coset of $H^1(C;\bbz_2)$
in $H^1(U(C);\bbz_2).$

Suppose $\omega$ is a smooth simple closed curve on $C.$ Let
$\overline{\omega}$ denote the lift of $\omega$ to $U(C)$ given by the unit
tangent vector at each point of $\omega.$ There are actually $2$ lifts, but
they are homotopic in $U(C)$ by the homotopy that rotates one tangential
direction to the other.

We will let $z\in H_1(U(C);\bbz_2)$ denote the class represented by the
tangential framing on $\partial(D^2),$ where $D^2$ is any closed $2$-disc in
$C$. As is nicely explained by Johnson, ``Intuitively a spin structure
$\zeta\in Spin (C )$ is a function assigning a number mod $2$ to each framed
curve of $C$, subject to the usual homological conditions and also that the
boundary of a disc in $C$ tangentially framed receives one".

\begin{lemma}(Johnson \cite{johnson}) Suppose the homology class $u\in
  H_1(C;\bbz_2)$ is represented by $\omega_1+\omega_2+\cdots + \omega_r,$
  where $\{\omega_1,\omega_2,\ldots,\omega_r\}$ is a set of non-intersecting
  smooth simple closed curves.  Then $
  \overline{\omega}_1+\overline{\omega}_2+\cdots + \overline{\omega}_r+rz$
  depends only on the class $u,$ and not on the particular representation.
\end{lemma}
\begin{definition}
$\tilde{u}= \overline{\omega}_1+\overline{\omega}_2+\cdots +
\overline{\omega}_r+rz.$
\end{definition}

This canonical lifting from $H_1(C;\bbz_2)$ to $H_1(U(C );\bbz_2)$ fails to be
a homomorphism. This is made explicit by the following result in
\cite{johnson}.

\begin{lemma}
  If $a,b\in H_1(C;\bbz_2)$ then $\widetilde{a+b}=\tilde{a}+\tilde{b}\
  +<a,b>z,$ where $<a,b>$ is the intersection pairing.
\end{lemma}

Let $e_1,\ldots,e_{2g}$ denote a basis of $H_1(C;\mathbb Z)\approx \mathbb
Z^{2g}$. We do not assume that this basis is symplectic. The following is
clear.

\ble A basis for $H_1(U(C);\mathbb Z)$ is given by
$\tilde{e}_1,\tilde{e}_2,\ldots, \tilde{e}_{2g},z.$ \ele

Let the dual basis of $H^1(U(C);\mathbb Z)$ be denoted
$\zeta_1,\ldots,\zeta_{2g},\eta.$ Then the mod 2 reduction $\bar{\eta}\in
H^1(U(C);\mathbb Z_2)$ is a particular spin structure on $C.$ It follows that
the set of spin structures is given by
$$\displaystyle
 Spin(C)=
\left\{\sum_{i=1}^{2g}x_i\bar{\zeta_i}+\bar{\eta}\ \bigm
|\ \mbox{all}\ x_i\in \mathbb Z_2\right\}.
$$
We can now determine the action of an automorphism $f:C\to C$ on
$Spin(C)$. Write
$f_*(e_i)=\sum_{j=1}^{2g}a_{ji}e_j,\ i=1,\ldots,2g,$
where the $a_{i,j}$ are the entries of $A\in SL_{2g}(\mathbb Z)$.

\begin{definition}
$v_i=\sum_{1\le j_1<j_2\le 2g}a_{j_1i}a_{j_2i}<e_{j_1}, e_{j_2}>,$
$V=[v_1,v_2,\ldots,v_{2g}].$
\end{definition}

We use the notation $V_f$ or $V_A$ if we want to emphasize that the vector $V$
comes from $f$ or $A.$ Note that if $f$ is a diffeomorphism of $C$ inducing
$f_*$ (resp. $f^*$) on $H_1$ (resp. $H^1$) then $\widetilde{f_*(a)} =
f_*(\widetilde{a})$ (and same for $f^*$). This is because $f_*(\overline a) =
\overline{f_*(a)}$ since $f_*$ acts by its differential on the tangent space
and because $f_*(z)=z$.  The following easy computation is an immediate
corollary of this fact

\begin{lemma}
With $f_*(e_i)=\sum_{j=1}^{2g}a_{ji}e_j$ for $i=1,\ldots,2g,$ we have
\begin{eqnarray*}
\displaystyle f_*(\tilde{e_i})&=&\sum_{j=1}^{2g}a_{ji}\tilde{e_j}+
v_iz\
\mbox{and}\ f_*(z)=z\\
f^*(\zeta_i)&=&\sum_{k=1}^{2g}a_{ik}\zeta_k\ \mbox{and}\
f^*(\eta)=\sum_{k=1}^{2g}v_k\zeta_k+\eta.
\end{eqnarray*}
\end{lemma}

Thus the matrices for $H_1(U(C);\mathbb Z)\stackrel{f_*}{\to} H_1(U(C);\mathbb
Z)$ and $H^1(U(C);\mathbb Z)\stackrel{f^*}{\to} H^1(U(C);\mathbb Z)$ are
$\left[
\begin{array}{cc}A& 0\\
V&1
\end{array}
\right]
$ and
$\left[
\begin{array}{cc}A^T& V^T\\
0&1
\end{array}
\right]
$ respectively.    A routine calculation gives

\begin{lemma}
 Suppose $\displaystyle \xi=\sum_{i=1}^{2g}x_i\bar{\zeta_i}+\bar{\eta}$ is a spin
structure. Then
$$f^*(\xi)=
\sum_{k=1}^{2g}\left(\sum_{i=1}^{2g}a_{ik}x_i+v_k\right)\bar{\zeta_k}+
\bar{\eta}.$$
\end{lemma}

Let $X$ denote the column vector $(x_1,\ldots,x_{2g})^T.$ Then

\begin{corollary}\label{invspinstrco}
$\displaystyle \xi=\sum_{i=1}^{2g}x_i\bar{\zeta_i}+\bar{\eta}$ is an invariant spin
structure for  $f:C\to C$ if, and only if,
\begin{eqnarray}\label{invspineq} (\bar{A}^T-I)\bar{X}= -\bar{V}^T
\end{eqnarray}
\end{corollary}

This is equivalent to $(A^T-I)X\equiv -V^T (mod\ 2).$ Proposition
\ref{mainth1} of \S \ref{intro} is now a consequence for if we suppose that
$\xi_1,\xi_2$ are invariant spin structures associated to column vectors
$X_1,X_2$, then each vector satisfies equation (\ref{invspinstrco}), and
therefore $X=X_1-X_2$ satisfies $(\bar{A}^T-I)\bar{X}^T=0.$

\subsection{Proof of Main Theorems \ref{mainth3} and \ref{mainth2}}
We first need information on the
similarity class of the matrix $A\in SL_{2g}(\mathbb Z).$

\begin{definition}
Let $\Phi_d(x)$ be the cyclotomic polynomial generated by the primitive $d^{th}$
 roots of unity and let $C_d$ denote the $\phi(d)\times \phi(d)$ companion matrix
 of $\Phi_d(x),$ that is
$$\displaystyle
C_d=\left[
\begin{array}{cccccc}0&1&0&0&\cdots&0\\
                     0&0&1&0&\cdots&0\\
                     \cdot&\cdot&\cdot&\cdot&\cdots&1\\
                     -a_0&-a_1&\cdot&\cdot&\cdots&-a_{\phi(d)-1}
\end{array}
\right]
$$
where $\phi(d)$ denotes Euler's totient function and
$\Phi_d(x)=a_0+a_1x+\cdots +a_{\phi(d)-1}x^{\phi(d)-1}+x^{\phi(d)}.$
\end{definition}

Suppose $f:C\to C$ is an automorphism of order $n.$ If $g\ge 2$ then $A$ also
has order $n.$ It follows that over the rationals $\mathbb Q,$ $A$ is similar
to a direct sum of companion matrices.  In fact, according to \cite{sjerve},
$\exists$ unique distinct divisors $1\le d_1<d_2\cdots<d_r\le n$ of $n$ and
unique positive integers $e_1,\ldots,e_r$ such that
\begin{itemize}
\item $\displaystyle A\ \mbox{is similar to}\ e_1C_{d_1}\oplus
  e_2C_{d_2}\oplus\cdots\oplus e_rC_{d_r}$ as a matrix in $GL_{2g}(\mathbb
  Q).$
\item $n=LCM(d_1,\cdots,d_r)$ and
$2g=e_1\phi(d_1)+e_2\phi(d_2)+\cdots e_r\phi(d_r).$
\end{itemize}

The minimal and characteristic polynomials of $C_{d}$ are both $\Phi_d(x).$ It
follows that the minimal and characteristic polynomials of $A$ are
$$\displaystyle
\mu_A(x)=\prod_{i=1}^{r}\Phi_{d_i}(x),\
\gamma_A(x)=\prod_{i=1}^{r}\Phi_{d_i}(x)^{e_i}\ \mbox{respectively.}
$$From the factorization
$\displaystyle x^{n-1}+x^{n-2}+\cdots +1=\prod_{d|n,d>1}\Phi_d(x) $ it follows
that $n=\prod_{d|n,d>1}\Phi_d(1)$ and
$det(I-A^T)=\gamma_A(1)=\prod_{i=1}^{r}\Phi_{d_i}^{e_i}(1).$ The following
lemma is well known.
\begin{lemma}
$\Phi_d(1)=\left\{
\begin{array}{ll}0&\ \mbox{if $d=1$}\\p&\ \mbox{if}\ d=p^k,\ \mbox{$p$ = a prime}\\
                 1&\ \mbox{in all other cases}
\end{array}\right.$
\end{lemma}

\begin{corollary}\label{det}
$\bar{A}^T-I$ is invertible (as a matrix over $\mathbb Z_2$) if, and only if,
$d_1>1$ and none of the $d_i$ are  powers of $2.$
\end{corollary}

From proposition (\ref{mainth1}) it follows that $f:C\to C$ has a unique
invariant spin structure if, and only if, $\bar{A}^T-I\in SL_{2g}(\mathbb
Z_2)$.  This happens if, and only if, $det(A-I)$ is odd. Now theorem
(\ref{mainth2}) of \S \ref{intro} will follow from this observation.

\bth  Suppose $f:C\to C$ is an automorphism of  order $n,$ where $n$ is odd. Then $f$ leaves only one spin structure invariant if, and only if, the associated orbit
  surface $C/\mathbb Z_n$ has genus zero.
\end{theorem}

\begin{proof}
If $V$ is a vector space and $L:V\to V$ is
a linear map let's write $E(L,\lambda)=\{v\in V \Bigm | L(v)=\lambda v\}$ 
the eigenspace of $L$ associated to $\lambda.$ Let
$\mathcal H$ (resp. $\mathcal H^*$) denote the complex vector space of holomorphic
(resp. anti-holomorphic) differentials on the curve $C$, $\dim{\mathcal H} = g$ 
the genus of $C$. The automorphism
$f:C\to C$ induces linear automorphisms $L:\mathcal H\to \mathcal H,$
$L^*:\mathcal H^*\to \mathcal H^*$. The subspaces
$E(L,1)$ and $E(L^*,1)$ correspond to the $\bbz_n$-invariant differentials,
where $\bbz_n$ is the cyclic subgroup generated by $f$, and thus
they have dimension $h:=genus(C/\mathbb Z_n).$ Now there is an
isomorphism $H^1(C;\mathbb Z)\otimes \mathbb C\approx \mathcal H\oplus
\mathcal H^*$ which is compatible with $f^*:H^1(C;\mathbb Z)\otimes \mathbb
C\to H^1(C;\mathbb Z)\otimes \mathbb C$ on the one hand, and with $L\oplus
L^*:\mathcal H\oplus \mathcal H^*\to \mathcal H\oplus \mathcal H^*$ on the
other.  From this it follows that the dimension of $E(f^*,1)=E(A,1)$ equals
$2h.$ From the rational canonical form $ e_1C_{d_1}\oplus
e_2C_{d_2}\oplus\cdots\oplus e_rC_{d_r}$ we see that the dimension of
$E(A,1)$ is $0$ if $d_1>1$ and $e_1$ otherwise.  The result now follows from
corollary \ref{det} and proposition \ref{mainth1}.
\end{proof}

Next we give a proof of theorem \ref{mainth3}.

\begin{proof} Suppose $C$ is hyperelliptic and $f$ is the hyperelliptic
  involution $J.$ Then the induced isomorphism $J_*:H_1(C;\mathbb Z)\to
  H_1(C;\mathbb Z)$ is $-I.$ Therefore $V_J=[0,0,\ldots,0]$ and equation
  (\ref{invspineq}) for an invariant spin structure becomes trivial.

  Conversely, suppose $f:C\to C$ is a non-trivial automorphism leaving every
  spin structure invariant. Then equation (\ref{invspineq}) is valid for all
  vectors $X$ and therefore $\bar{A}= I.$ By a theorem of Serre (see p. 275 of
  \cite{farkas}) we see that $A$ must have order $2,$ and therefore so does
  the automorphism $f:C\to C.$ Now a theorem of Nielsen \cite{nielsen} states
  that $f$ is determined up to conjugacy by its fixed point data.  In the case
  of an automorphism of order $2$ the fixed point data is just the number of
  fixed points, and therefore two involutions with the same number of fixed
  points are conjugate.  Let the number of fixed points be $r.$ It is known
  \cite{farkas} that if $f\neq id$, $r\leq 2g+2$.

\begin{figure}[htb]\label{spinfig2}
\begin{center}
\epsfig{file=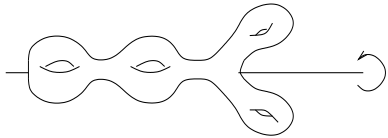,height=0.7in,width=2in,angle=0}
\caption{An involution with $r<2g+2$ fixed points}
\end{center}
\end{figure}

Assume $r<2g+2.$ Let $2s=2g+2-r>0$ ($r$ must be even because of the
Riemann-Hurwitz formula).  Then $f$ must be conjugate to an automorphism as
depicted in the Figure 1, and therefore we can choose a basis so that
$$\dis
A\ =\left[\begin{array}{ccc}
-I_r&0&0\\
0&0&-I_s\\
0&-I_s&0
\end{array}\right]
$$
But this contradicts the equation $\bar{A}= I,$ and therefore $r=2g+2.$
Thus $C$ is hyperelliptic and $f$ is the hyperelliptic involution.
\end{proof}

\subsection{Corollaries}
Corollary \ref{invspinstrco} and the proof of theorem \ref{mainth3} have some
elementary consequences. Let $\rho : Aut(C)\lrar Sp_{2g}(\bbz_2)$ be the
homomorphism induced from the action on $H_1(C;\bbz_2)$ where $Sp_{2g}(\bbz)$
is the symplectic modular group and $Sp_{2g}(\bbz_2)$ its mod-$2$ reduction.

\bco\cite{farkas} $\rho$ is a monomorphism if and only if $C$ is
non-hyperelliptic.  In the hyperelliptic case $ker(\rho )=\bbz_2$ is generated
by the hyperelliptic involution.  \eco

Let $\Gamma_g(C)$ be the mapping class group of $C$; that is the group of
isotopy classes of orientation preserving diffeomorphisms of $C$.  It is
well-known that $\Gamma_g(C)$ is a finitely generated group with generators
the Dehn twists around representative loops of a symplectic basis of $C$. If
$[f]\in\Gamma_g$, then the action of $f$ on $Spin(C)$ (viewed as cohomological
classes) is independent of the choice of this representative so that
$\Gamma_g$ acts on $Spin(C)$ .  Let $S_g\subset\Gamma_g$ be the subgroup of
$\Gamma_g$ that fixes all spin structures. The following is immediate

\bco\cite{sipe} The subgroup $S_g$ is precisely the subgroup of elements
that induce the identity on $H_1(C;\bbz_2)$.
\eco

In fact in \cite{sipe}, Sipe states her theorem for all $n$-roots and she
identifies $S_g$ with the subgroup of elements which induce the identity on
$H_*(U(C);\bbz_2)$.  But an orientation preserving automorphism is the
identity on $H_*(U(C);\bbz_2)$ if and only if it is the identify on
$H_*(C;\bbz_2)$.

%******************************************************************************

\section{Spin Structures on Hyperelliptic Curves}\label{hyperelliptics}

We use a convenient description of the spin structures on a hyperelliptic
curve in terms of divisors due to Mumford (\cite{mumford} or \cite{arbarello},
appendix B) to count invariant spin structures for subgroups of the 
automorphism group.

From the appendix we have that
\begin{equation}\label{description}
Spin(C ) = \{D\in Cl(C )\ |\ 2D = K\}
\end{equation}
where $Cl(C )$ is the divisor class group of $C$.
We will denote by 
$J_2(C )$ the subgroup of points of order two in $J(C )=Cl^0(C ),$
the ``Jacobian" of line bundles of degree zero. See \cite{arbarello}.

\subsection{Spin Divisors on a Hyperelliptic Curve}\label{spinonhyp}
We consider the hyperelliptic surface $y^2 =
\prod_{i=1}^{2g+2}(x-e_i)$ having genus $g$ and
branch set $B=\{e_1,\ldots,e_{2g+2}\}$.  The $e_i$ are distinct points
in the Riemann sphere. We will write $\pi : C\lrar\bbp^1$ for the degree
two covering map sending $(x,y)\mapsto x$, and we let $p_i$ be the ramification points
such that $\pi (p_i)=e_i$.
A description of the canonical divisor $K$ of $C$ is completely standard.
Let $D$ be any divisor of the form $2p_i$ or $x+y$ if
$x,y\not\in\{p_1,\ldots, p_{2g+2}\}$ (i.e. $D$ is the divisor associated to the
line bundle $L$ obtained as the pullback of the unique line bundle over
$\bbp^1$ of degree $+1$). Then \cite{arbarello, mumford}
\begin{equation}\label{canonical} 
K_C = -2D + p_1+\cdots +p_{2g+2} = (g-1)D
\end{equation}

We can next determine the divisor of an element $L\in J_2(C )$.  Let
$\phi : \bbp^1\lrar\bbc$ be the meromorphic function $\phi (z) = (z-e_i)/
(z-e_j), i\neq j$, and consider the composite
$C\fract{\pi}{\lrar}\bbp^1\fract{\phi}{\lrar}\bbc$. The divisor of the
composite is $(\phi\circ\pi) = 2p_i-2p_j$ and hence by definition $p_i-p_j$ is
the divisor of a line bundle of order two on $C$ ($i\neq j$). Using the
same argument with the rational function $\phi (z ) = {(z-e_i)(z-e_j)/ z^2}$
shows that $E = p_i-p_j= D-p_i-p_j$. Generally one sees that if $L\in
J_2(C )$, then it is represented by either one of the divisors
\begin{equation}\label{representation}
E= p_{i_1}+ \cdots + p_{i_k}- p_{j_1}-\cdots - p_{j_k}=
kD - p_{l_1}-\cdots -\cdots -p_{l_{2k}}
\end{equation}
Note that the above formula shows that any finite set $T\subset\{1,\ldots,
2g+2\}$ of even cardinality gives an element in $J_2(C)$. One can further see that any such
$T$ and its complement $T^c\subset B$ give rise to the same element. 
If we define $E_g$ to be the quotient 
$$E_g:=\{T\bigm | T\subseteq B, |T|\ \hbox{even} \}/\sim,\ \mbox{where}\
T^{\prime}\sim T\Longleftrightarrow T^{\prime}=T\ \mbox{or}\ T^{\prime}=T^c$$
Then the map
$$\alpha : E_g\lrar J_2(C)\ \ ,\ \ T\mapsto \alpha_T = \sum_{i\in T}(p_i - p_{2g+2})$$
defines an automorphism between two copies of $\bbz_2^{2g}$ \cite{dolgachev}.

\bth\label{spinonhyper}\cite{mumford} Let $C$ be as above, $p_i\in C$. Then
every theta characteristic is of the form
$$\theta_T = kD + p_{i_1} + \cdots + p_{i_{g-1-2k}}$$
for $-1\leq k\leq {g-1\over 2}$,
$T := \{p_{i_1}, \cdots , p_{i_{g-1-2k}}\}$ with the $i_{\alpha}$ distinct.
Moreover such a representation is unique if $k\geq 0$ and subject to a single
relation $-D+ p_{i_1} + \cdots + p_{i_{g+1}}=
-D+ p_{j_1} + \cdots + p_{j_{g+1}}$
when $k=-1$.
\end{theorem}

\begin{proof} An easy sketch is in \cite{arbarello}, exercises 26-32, p.
  287-288.  We work in the divisor class group where we have the (equivalent)
  relations $p_{i_1} + \cdots + p_{i_{g+1}}= p_{j_1} + \cdots + p_{j_{g+1}}$
  and $K= (g-1)D$.  With $E$ as in the theorem we can then compute
\begin{equation}\label{relation}
2E= 2kD + 2p_{i_1} + \cdots + 2p_{i_{g-1-2k}}=
2kD + (g-1-2k)D = (g-1)D= K
\end{equation}
which implies that $E$ is a spin divisor.
The previous paragraphs also show that
$$p_{i_1}+\cdots +p_{i_r}=  p_{j_1}+\cdots +p_{j_r}\
\Longleftrightarrow\ r=g+1$$
for distinct $i_{\alpha}$'s and $j_{\beta}$'s
so that a representation of $E$ as in (\ref{relation}) is unique
when $0\leq k\leq (g-1)/2,$ thus yielding $\sum_{k=0}^{(g-1)/2}{2g+2\choose
g-1-2k}$ distinct spin structures. The one relation for $r=g+1$ yields 
${1\over 2}{2g+2\choose
g+1}$ additional structures. Since this adds up to exactly $2^{2g}$, we have
accounted for all spin structures this way.
\end{proof}

Since to each branched point corresponds a unique ramification point, we can
write

\ble\label{description2}
$Spin(C)$ is the set of all $T\subset B$ such that $|T|\equiv g+1\mod (2)$
modulo the equivalence relation $T\sim T^c$. 
This has a natural affine structure over $E_g$ given by
$\theta_S + \alpha_T = \theta_{T+S}$.
\ele

%**************************************************************************

\subsection{Automorphisms}\label{auts}
Let $f:C\lrar C$ be an automorphism of a Riemann surface.
Then $f$ acts on $Spin (C )$ as in (\ref{description}) by sending
$$f(\sum n_iP_i) = \sum n_if^{-1}(P_i)$$
Indeed one way to see this is to replace the cotangent bundle by the tangent bundle,
and since $f$ is an automorphism, replace the pullback $f^* : T^*_{f(x)}C\lrar
T^*_xC$ by $f^{-1}: T_{f(x)}C\lrar T_xC$.

For the rest of this section $C$ will be hyperelliptic, $B$ its branched set
as in section \ref{spinonhyp} and $J$ the hyperelliptic involution. This is a
central element in $Aut(C)$ and we have a short exact sequence
$$0\lrar\bbz_2\{J\}\lrar Aut(C)\lrar\overline{Aut}(C)\lrar 0$$
where $\overline{Aut}(C)$ is a finite subgroup of $PSL_2(\mathbb C)$ and
hence is a cyclic, dihedral or  polyhedral group.
If $f:C\to C$ is an automorphism then we associate
$\bar f\in \overline{Aut}(C)$. Here $f$ acts
on the ramification (i.e. Weierstrass) points and $\bar f$ on the branch set $B$.
Since an element of $PSL_2(\mathbb C)$ other than the identity fixes at most 
two points, we see that if
$f\neq id,J$,
then $\bar f$ acts on $B=\{e_1,\ldots, e_{2g+2}\}$ with at most two fixed points. 

We will be then distinguishing \underline{three cases}: when $\bar f$ acts fixed point
freely, with one fixed point or with two.

We assume below that $f\neq id,J$, $f$ has order $n$ so that (wlog) $\bar f
(e) = \zeta e$, $\zeta = e^{2\pi\imath/n}$.  We use the notation $<e>:=\{e,\zeta
e,\ldots, \zeta^{n-1}e\}$ to denote the orbit of $e\in\mathbb P^1(\mathbb
C)=\bbc\cup\{\infty\}$ under the action of the cyclic group $\mathbb
Z_n\subset PSL_2(\mathbb C)$ generated by $\bar f$.  Every orbit $<e>$ has length $n$
except possibly for two which are singletons: $<0>=\{0\}$ and
$<\infty>=\{\infty\}.$

A spin structure $[T]$ as in lemma \ref{description2} is invariant under the
automorphism $f:C\to C$ if, and only if, $\zeta T=T\ \mbox{or}\ T^c.$ We now put
this to good use.

\begin{proposition}\label{hyperinvodd}
Assume $f\neq id$ has odd order $n$.  
Then the number of invariant spin structures under $f$ is
$$
\left\{\begin{array}{ll}
2^{r-2}&if \ \bar f\ \mbox{acts freely},\ \ \mbox{where}\ 2g+2=nr \\
2^{r-1}&if \ \bar f\ \mbox{acts with one fixed point},\ \ \mbox{where}\ 2g+2=nr+1\\
2^{r}&if \ \bar f\ \mbox{acts with two fixed points},\ \ \mbox{where}\ 2g=nr\end{array}\right.
$$
\end{proposition}

\begin{proof}
  Since $n$ is odd an invariant spin structure is determined by a subset
  $T\subset B$ such that $\zeta T=T$ and $card(T)\equiv\ g+1\ (mod\ 2).$
  Assume $\bar f$ acts freely on $B$ so that $B$ is the disjoint union
$$B = <e_1>\sqcup <e_2>\sqcup\cdots\sqcup <e_r>,$$ 
where all the $e_i$ are in $\bbc^*=\bbc-\{0\}$.  By this notation we also mean
that we have reordered the $e_i$ so that $e_1,\ldots,e_r$ are representatives
of the orbits $<e_i>$ for $i=1,\ldots,2g+2.$ There is no loss of generality in
doing this.  In that case and if some $\zeta^m e_i\in T$, then the entire
orbit $<e_i>\subset T.$ Therefore $T$ must be a disjoint union of the form
$$\dis
T=<e_{i_1}>\sqcup<e_{i_2}>\sqcup\cdots\sqcup <e_{i_k}>,\ \mbox{where
$1\le i_1<i_2<\cdots < i_k\le r$ and $k= 0,\ldots,r$}.
$$

Therefore $card(T)=kn\equiv k\ (mod\ 2).$ The number of such $T$ is
$\dis\sum_{k=0}^{r} {r\choose k}=2^r,$ and the number satisfying
$card(T)\equiv g+1\ ( mod\ 2)$ is $2^{r-1}.$ This follows from the identity
$\dis\sum_{k=0}^{r}(-1)^k {r\choose k}=0$. These subsets come in complementary
pairs $\{T,T^c\},$ and therefore the number of invariant spin structures is
$2^{r-2}$.

Assume $\bar f$ acts with one fixed point so that
$$B = <\cdot>\sqcup <e_1>\sqcup <e_2>\sqcup\cdots\sqcup <e_r>,\ \ 
e_i\in \bbc^*, \forall i$$ 
where $<\cdot>$ is either $<0>$ or $<\infty >$.
If $T\subset B$ satisfies $\zeta T=T,$ then the
 possibilities for $T$ are
\begin{eqnarray*} T&=&\dis <e_{i_1}>\sqcup<e_{i_2}>\sqcup\cdots\sqcup <e_{i_k}>,\
\mbox{or}\\
T&=&\dis <\cdot >\sqcup<e_{i_1}>\sqcup<e_{i_2}>\sqcup\cdots\sqcup <e_{i_k}>
\end{eqnarray*}
In both cases we must have  $1\le i_1<i_2<\cdots <i_k\le r$ and $k=0,\ldots,r.$
The number of such subsets $T$ is $2^{r+1},$ and the number satisfying $card(T)\equiv
g+1\ (mod\ 2)$ is $2^r.$  These subsets come in complementary pairs, and therefore the
number of invariant spin structures is $2^{r-1}.$ 

Finally assume that $\bar f$ has two fixed points in $B$ so that
$$B = <0>\sqcup<\infty>\sqcup <e_1>\sqcup <e_2>\sqcup\cdots\sqcup <e_r>,\ 
e_i\in \bbc^*, \forall i$$ 
If $T\subset B$ satisfies $\zeta T=T$ then $T$ can be any collection of orbits
with the right cardinality so that by arguments similar to those above we also
get $2^r$ invariant spin structures.
\end{proof}

\begin{proposition}\label{hyperinveven}
  Assume $f\neq id, J$ has even order $n$. Then the number of invariant spin structures is
$$\left\{\begin{array}{ll}
 2^{r-1}&if\ \bar f\ \mbox{acts freely on $B$ and $g$ is even}\ \ (2g+2=nr)\\
                          
 2^r&if\ \bar f\ \mbox{acts freely on $B$ and  $g$ is odd}\ \ (2g+2=nr)\\
                        
 2^r&if\ \bar f\ \mbox{acts on $B$ with two fixed points}\ \ (2g=nr)
\end{array}\right.
$$
\end{proposition}

Note that when $n$ is even, $\bar f$ cannot act on $B$ with 
a single fixed point since $2g+2=nr+1$.

\begin{proof}
  First assume $\bar f$ acts on $B$ fixed point freely so that again $B$ is
  the disjoint union $<e_1>\sqcup <e_2>\sqcup\cdots\sqcup <e_r>,$ where all
  the $e_i$ are in $\bbc^*=\bbc-\{0\}$.  Then there are 2 possibilities for an
  invariant spin structure $[T]:$ either $\zeta T=T^c$ or $\zeta T=T.$

  If $\zeta T=T^c$ occurs then necessarily $T$ is the disjoint union of the
  sets
$$\{\zeta^{\epsilon_j} e_j,\zeta^{2+\epsilon_j} e_j,\zeta^{4+\epsilon_j} e_j,\ldots,
\zeta^{n-1+\epsilon_j}e_j\},\ \mbox{where $j=1,\ldots,r$ and  $\epsilon_j=\pm 1$
$\forall j$}$$
This means that $n$ is necessarily even (as is the case) and that
$\zeta$ fixes no branch points.
The number of such subsets is $2^r$ since it equals the number of choices of the
$\epsilon_j.$ Every such subset $T$ determines a spin structure $[T]$ because
$card(T)=g+1$ in this case. These subsets come in equivalent pairs $\{T,T^c\},$
and therefore there are $2^{r-1}$ invariant spin structures in this case.

For the second possibility there are $2^r$ subsets $T\in \mathcal S_B$ such
that $\zeta T=T,$ namely
$$T=\dis <e_{i_1}>\sqcup<e_{i_2}>\sqcup\cdots\sqcup <e_{i_k}>,\
\mbox{where $1\le i_1 <\cdots < i_k\le r,$ $k=0,1,\ldots,r$}$$
Now $card(T)=kn\equiv 0\ (mod\ 2),$ and therefore $[T]$ is a spin structure if, and
only if, $g$ is odd.  These spin structures come in complementary pairs, so we have
$2^{r-1}$ invariant spin structures for this possibility.

In summary, when $\bar f$ acts fixed point freely 
there are $2^{r-1}$ invariant spin structures if $g$ is
even, and  $2^{r-1}+2^{r-1}=2^r$ if $g$ is odd.

The arguments for the other case are similar.
\end{proof}

%******************************************************************************

\subsection{Genus Two Curves}\label{genus2}

As an exercise we count the invariant spin structures for all possible
subgroups of $Aut(C)$ in the case when $C$ is of genus two; i.e.
$C : y^2 =f(x)$ where $f$ is a polynomial of
degree $5$ or $6$. A complete list of all possible reduced automorphism
groups that can occur together with their associated equations
has been long known (Bolza). For those groups we can count precisely the
number of fixed spin structures. It is given by the table below.

\begin{table}[ht]
%\caption{}\label{eqtable}
\renewcommand\arraystretch{1.5}
\noindent\[
\begin{array}{|c|c|c|c|}
\hline
g=2&\overline{Aut}&\mbox{representative}\ f&\fract{\hbox{\# of fixed spin}}
{\mbox{structures}}\\
\hline
(i)&\bbz_2&(x^2-1)(x^2-a)(x^2-b)&4\\
\hline
(ii)&D_4&(x^2-1)(x^2-a)(x^2-{1\over a})&2\\
\hline
(iii)&D_6&x^6-2ax^3+1&1\\
\hline
(iv)&D_{12}&x^6+1&1\\
\hline
(v)&S_4&x(x^4-1)&0\\
\hline
(vi)&\bbz_5&x^5-1&1\\
\hline
\end{array}
\]
\end{table}

Here the actual new statement is in the last column.  
It is easy to see that the groups in the table are the only ones allowed since
besides being finite in $SO(3)$,
they also are subgroups of $S_6$ (permutations of the six branch points).
Note that the spin
structures in genus 2 are given by all equivalence classes of subsets of cardinality 3 in $B$ ($10$
of them and these are the even structures) together with all subsets of
cardinality $5$.  For each case in the table above, the configuration of
points of $B$
in $\bbp^1$ is laid out according to the geometry of the associated
group action. It is easier to write $B=\{1,2,3,4,5,6\}$.

\begin{proof} (last column of table)\\
$\bullet$ (i) $\bbz_2$ acts fixed point freely via the map $x\mapsto -x$ and
hence according to proposition \ref{hyperinveven} it fixes $2^{r-1}=4$ branch
points (here $n=2,r=3$).\\
$\bullet$ (ii) This has to do with the geometry of a parallelogram with vertices at
$a,\bar a,-\bar a,-a$ say ($|a|=1$). Let's refer to these branch points by $1,2,3,4$
respectively (so that $``5"$ and $``6"$ correspond to the remaining branch points 
$1$ and $-1$). The two involutions are $\tau_1 = (12)(34)$ (reflection through the
$x$-axis) and $\tau_2 = (13)(24)(56)$ (reflection through the $y$-axis). 
According to proposition \ref{hyperinveven},
$\tau_1$ and $\tau_2$ fix $4$ spin structures each but it is easy to see that
they have only two such structures in common namely $\{5,3,4\}$ and
$\{5,1,2\}$.
\\
$\bullet$ (iii) 
This is the situation where the branch points are configured in the plane as
in the figure below (each dot represents a branch point in $\bbc\subset\bbp^1$).

\begin{figure}[htb]
\begin{center}
\epsfig{file=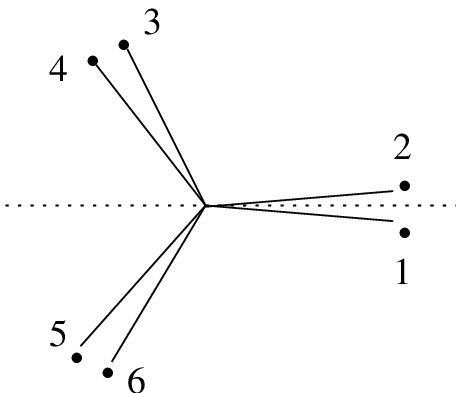,height=1.3in,width=1.5in,angle=0.0}
%\caption{}
\end{center}
\end{figure}

The polynomial in (iii) has roots $e^{\pm\alpha}, je^{\pm\alpha},
j^2e^{\pm\alpha}$ numbered as in the figure. There is only one spin structure
invariant by the $2\pi\over 3$ rotation: namely $T=\{1,3,5\}\sim
T'=\{2,4,6\}$. The involution with respect to the $x$-axis switches $T$ and
$T'$ and
hence $D_6$ has a unique invariant spin structure.\\
$\bullet$ (iv) $\bbz_6$ acts cyclically on $B$ without fixed points and hence
according to proposition \ref{hyperinveven} it has a single fixed spin
structure
$T=\{1,3,5\}\sim T^c$.\\
$\bullet$ (v) the branch points form an octahedron $\{1,2,3, 4, N,S\}$ with
north $N$ and south $S$ poles.  The cylic $\bbz_4$ acts by rotating the
equator $1\mapsto 2\mapsto 3\mapsto 4$.  According to proposition
\ref{hyperinveven} this group fixes $2$ spin structures which must be odd of
the form $\{1,2,3, 4, N\}$ or $\{1,2,3, 4, S\}$.  But there is an involution
switching $N$ and $S$ and thus
there can't be any invariant spin structure.\\
$\bullet$ (vi) corresponds to when $\overline{Aut}$ is of order $n=5$
generated by $\bar f$ (in the notation of section \ref{auts}) which acts by
fixing one point $\infty$. According to proposition \ref{hyperinvodd}
$\bar f$ fixes $2^{r-1}=1$ structures.
 \end{proof}

%******************************************************************************

\section{Spin Structures on Klein's Curve}\label{kleincurve}

In this section we show that the Klein's quartic curve $\mathcal K$
has a unique spin structure invariant under all automorphisms.
Let $G=PSL(2,{\bf F}_7)=Aut({\mathcal K}).$ This is a simple group with presentation of the form
$$\displaystyle
G=\left<R,S,T \bigm | R^2=S^3=T^7=RST=1,\ \mbox{etc}\right>
$$

Let $e_1,e_2,\ldots,e_6$ be a standard symplectic basis of $H^1(C;\mathbb Z).$ That
is a basis of the free abelian group  $H^1(C;\mathbb Z)$ chosen so that the
intersection pairing is given as follows:

$$ <e_i,e_j>=
\left\{
\begin{array}{lll}
+1 & \mbox{if $j=i+3$}\\
-1 & \mbox{if $j=i-3$}\\
~0 & \mbox{in all other cases}
\end{array}
\right.$$

From the work of Rauch and Lewittes 
\cite{rauch} we can find matrix representations for the induced maps
$R_*,S_*,T_*$ on $H_1(C;\mathbb Z):$
$$
R_*=\left[\begin{array}{cccccc}
-1&0&0&0&0&0\\
-1&1&-1&0&0&0\\
0&0&-1&0&0&0\\
0&0&0&-1&-1&0\\
0&0&0&0&1&0\\
0&0&0&0&-1&-1\\
\end{array}\right]$$
$$\ S_*=\left[\begin{array}{cccccc}
-1&0&0&-1&-1&-1\\
-1&0&1&0&0&-1\\
-1&1&0&0&0&0\\
1&0&0&0&0&0\\
0&-1&1&0&0&0\\
0&1&-1&0&1&1
\end{array}\right]$$
$$\ T_*=\left[\begin{array}{cccccc}
0&0&0&-1&-1&0\\
0&0&-1&-1&-1&0\\
0&0&-1&-1&0&0\\
1&0&0&1&1&1\\
-1&1&0&0&-1&-1\\
1&-1&0&0&1&0
\end{array}\right]
$$

The corresponding $V$ vectors are
$$V_R=[0,0,0,0,0,0],V_S=[-1,1,1,0,0,0],V_T=[0,0,0,-1,0,0].
$$
Using equation (\ref{invspineq}), together with a Maple program,  it is now
direct to determine the sets of spin
structures invariant under $R,S,T$ respectively. The solution vectors $\bar{X}$ are:
\begin{enumerate}
\item $\bar{X}=[x_1,0,x_3,x_4,x_5,x_4],$ where the $x_j$ are arbitrary, for $R$-invariant
spin structures.
\item $\bar{X}=[x_1+x_2,x_1+x_2,1,x_1+x_2,x_1,x_2],$ where the $x_j$ are arbitrary,
for $S$-invariant spin structures.
\item $\bar{X}=[0,0,1,0,0,0]$ for $T$-invariant spin structures.
\end{enumerate}

It is now clear that Klein's curve admits only one invariant spin
structure.  This is the spin structure fixed by $T$ and which by choosing the
$x_j$ appropriately we can show is also an invariant for $R$ and $S$.

\subsection{Spin Structures on Quartics}
We would like to say more about what this particular spin structure on
${\mathcal K}$ is.  Generally let $C$ be a quartic curve in $\bbp^2$.  There
are $2^6= 64$ spin structures on $C$ as pointed out, of which $28$ are ``odd''
and $36$ ``even'' (see appendix).  The canonical divisor class on $C$ is the
locus of the intersection of a hyperplane $\bbp^1\subset\bbp^2$ with $C$.
Suppose $C$ has a bitangent going through the points $p,q$, $p\neq q$.  This
is a line $\bbp^1$ intersecting the curve in $p,q$ with multiplicity two, and
hence $K=2p+2q$ (in the divisor class group). By definition, $p+q$ is a spin
structure on $C$. Since there are $28$ bitangents, this accounts for all odd
spin structures.

The even spin structures on the other hand are described as follows.
Following \cite{mukai}, a secant line $\overline{ab}$, $a\neq b\in C$ is
a \textit{strict tangent} of $C$ if either:\\
(i) $\overline{ab}$ is tangent to $C$ at a point $r$ different from $a,b$,
or\\
(ii) $\overline{ab}$ is a triple tangent at $a$ or at $b$.

This is equivalent to the condition that the divisor class $K=a+b+2r$ for
a point $r\in C$. A triangle is now \textit{strictly biscribed} if all sides
are strict tangents (see figure 2). A plane cubic has no such triangles but
a plane quartic has $288$ of them! The way they arise as it turns out
is in groups of $8$ triangles, one for every even spin structure.

\begin{figure}[htb]
\begin{center}
\epsfig{file=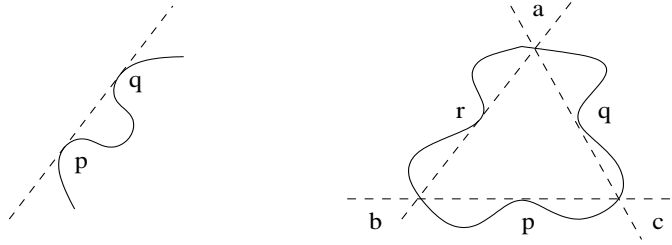,height=1.3in,width=3.5in,angle=0.0}
\caption{A bitangent (odd theta) and a strictly biscribed triangle (even theta)}
\end{center}
\end{figure}

Let $\Delta(a,b,c)$ be a strictly biscribed triangle in a quartic $C$ with
double tangents at $p,q,r$ as in figure 2. Then we can associate to $\Delta$
the even spin class
$$\theta (\Delta) = a+b+c+p+q+r - K$$
Using the relations $K=2p+b+c = 2q+a+c= 2r+a+b$ we can indeed check that
$2\theta (\Delta ) = K$.
It is known that for an even theta charateristic $\theta$, there are $8$
strictly biscribed triangles $\Delta$ such that $\theta(\Delta )=\Delta$
\cite{dolgachev, mukai}.

Restricting to the case of the Klein curve ${\mathcal K}$, we know that
$Aut({\mathcal K}) = PSL(2,{\bf F}_7)$ fixes no bitangent and hence an invariant
spin candidate must be even. Consider
the points $a=[1:0:0], b=[0:1:0]$ and $c=[0:0:1]\in {\mathcal K}$, $\mathcal
K$ being again the curve $x^3y+y^3z+z^3x=0$ in $\bbp^2$.
As discussed in \cite{bavard} for instance, the
transformation $T$ of order $7$ can be represented by
$$T : x\mapsto \gamma x\ \ ,\ \ y\mapsto \gamma^4y\ \ ,\ \ z\mapsto \gamma^2z$$
where $\gamma = e^{2\pi i\over 7}$. This transformation fixes the three
points above and hence the triangle $xyz=0$ they define. Note that
$a,b,c$ are inflexion points (i.e. points with triple tangents)
and hence they form what Klein calls an inflexion triangle; that is the
tangent at $a$ meets ${\mathcal K}$ at $b$, and the tangent at $b$ meets
${\mathcal K}$ at $c$ and finally the tangent at $c$ cuts ${\mathcal K}$ in
$a$.  This is also a biscribed triangle according to our definition with
associated theta divisor
$$\theta  := a+b+c+a+b+c - K = 2a+2b+2c - K$$
(since $r=a, p=b, q=c$). 

Let's check that $\theta$ is fixed by the generators $S$ and $R$ of order two
and three respectively. Now 
$$R: x\mapsto y \ \ ,\ \ y\mapsto z\ \ ,\ \ z\mapsto x$$
so $\Delta$ and hence $\theta$ is invariant. To see that $S$ acts trivially
on $\theta$ it doesn't quite help to write the long expression for $S$ given
also in \cite{bavard}. Instead consider the quotient $H={\mathcal K}/<S>$ by
the cylic subgroup generated by $S$ and write
$\pi : {\mathcal K}\lrar H$ the quotient map.
Clearly $H$ cannot be $\bbp^1$ since
$\mathcal K$ is not hyperelliptic. The only other possibility by the
Riemann-Hurwitz formula is that $genus(H)=1$ and $H$ is an elliptic curve. If
$a\in{\mathcal K}$, then $\pi^{-1}(\pi (a)) = \{a, Sa\}$. But on $H$ there
is a meromorphic function $f$ with a double pole at $\pi (a)$ and hence
in the divisor class group and by definition
$0 = (f\circ\pi ) = 2a-2S(a)$.
This means that $2a+2b+2c = 2S(a)+2S(b)+2S(c)$ and hence $S(\theta ) = \theta$
as desired. 

The above discussion proves theorem \ref{mainth4}.

\bre
In \cite{dolgachev2}, the action of $PSL(2,{\bf F}_7)$ on the set of even spin
structures is described with orbit decomposition $36= 1+7+7+21$. The one fixed
spin generator is expressed explicitly in terms of the orbital divisors as
$D_2-D_3-7D_7$, where $Aut(C)$ acts on ${\mathcal K}$ with orbits $D_2,D_3,D_7$ with
stabilizers of orders $2, 3, 7$.
\ere

%******************************************************************************

\section{Appendix : Cohomology and Divisors}\label{appendix}

In this appendix we reconcile both definitions of $Spin(C)$ given in
\S\ref{spin} and in \S\ref{hyperelliptics}.  Let $Cl(C)$ be the divisor class
group of $C$; that is the set of all linearly equivalent divisors. The $0$
divisor in $Cl(C)$ is the divisor of any meromorphic function on $C$. The
canonical divisor $K$ is the divisor associated to any meromorphic one form
$\omega $ on $C$. We write $K=(\omega )$. A spin structure is then a divisor
$D\in Cl(C)$ such that $2D = K$. We would like to see that this gives us a
cohomology class in $H^1(U(C);\bbz_2)$ which restricts to the generator in the
fiber.

The form $\omega$ is dual to a vector field $\mathcal V$ and the
``singularities'' of $\omega$ and $\mathcal V$ are the support of $D$; that is
the set of zeros and poles of $\omega$.  Note that a zero of $\omega$ is a
pole of $\mathcal V$ and vice-versa.  If $(\omega) = 2D$ then the index of the
vector field around any of its singular points is even.

Let $c$ be a closed immersed and framed curve in $C$ (this is precisely the
one cycles in $H_1(U(C ))$) and let us assume that $c$ avoids the singular set
of $\mathcal V$. We can then define the relative winding number of $\mathcal
V$ along the curve $c$ with respect to the framing of $c$.  This is in words
the number of times the vector field rotates in relation to this framing. More
specifically, set $v$ to be the framing, and choose a vector field $w$ such
that $(v,w)$ is a positively oriented basis along $c$.  
We write ${\mathcal V} = p_1v + p_2w$ where $p =
(p_1,p_2):S^1\lrar\bbr^2$.  We then define this winding number to be $\deg
p/|p|$ and we define $Ind_{\mathcal V}(c)\in\bbz_2$ to be this number modulo
two. This is a topological invariant and does not change under continuous
deformations of the vector field.

Observe at this point that if $c$ is a $1$-boundary (i.e. it bounds a framed
disk), then $Ind_{\mathcal V}(c) = 0$ precisely because the singularities of
${\mathcal V}$ have even index. On the other hand and if $c$ is a small simple
closed curve with tangential framing, then this framing cannot extend to the
interior disk bounded by $c$, and consequently the relative winding number is
odd and $Ind_{\mathcal V}(c) = 1$. Note also that a small closed tangentially
framed curve represents the generator of $H_1(fiber )$ in $H_1(U(C );\bbz_2)$.
The function $Ind_{\mathcal V}$ is then a well-defined cohomology class
$$Ind_{\mathcal V}\in H^1(U(C ),\bbz_2)=Hom(H_1(U(C
);\bbz_2),\bbz_2)=\bbz_2^{2g+1}$$ which restricts to the generator of the
cohomology of the fibers. This is a spin structure.

If $c$ is an immersed closed curve in $C$ that avoids the singularities
of $\omega$, then the associated quadratic form is given according to \cite{johnson}
$$q_{\omega}(c) = Ind_{\mathcal V}(c) + \# \mbox{of components of $c$} + \#\mbox{of
  self-intersections of $c$}$$

A spin structure $D\in Cl(C)$ is even (or odd) if the dimension of the vector
space of sections of $L(D)$ the line bundle associated to $D$ is even (or odd).

{\sc Acknowledgement:} The first author is grateful to Jean d'Almeida for 
interesting conversations regarding \cite{dolgachev}.

%*******************************************************************

\bigskip
\hfill\
{\footnotesize
\parbox{5cm}{Sadok Kallel\\
Laboratoire Painlev\'e,\\
Universit\'{e} de Lille I,
France\\
{\sf Sadok.Kallel@math.univ-lille1.fr}}\
{\hfill}\
\parbox{5cm}{Denis Sjerve\\
Department of Mathematics,\\
University of British Columbia,\\ Canada\\
{\sf sjer@math.ubc.ca}}}\\
\hfill

\end{document}